Cover letter

This manuscript is motivated by presenting a new combined numerical scheme based on backward differentiation formula (BDF) and generalized differential quadrature method (GDQM). The proposed scheme has been applied to the well-known KdV equation in different models. The stability analysis has been investigated, too.

# Efficient Numerical Scheme for Solving (1+1), (2+1)-dimensional and Coupled Burgers' Equation

N. A. Mohamed, A.S. Rashed[*]

*Department of Physics and Engineering Mathematics, Faculty of Engineering, Zagazig University, Egypt.*

**ABSTRACT**

A numerical scheme based on backward differentiation formula (BDF) and generalized differential quadrature method (GDQM) has been developed. The proposed scheme has been employed to investigate three cases of Burgers' equation, one-dimensional, two-dimensional and two-dimensional coupled models. The results showed an effectiveness accuracy in absolute error and error norms $L_2$ and $L_\infty$ compared to other methods. The results were evaluated at different values of Reynold's number, *Re*, and viscosity, $v$.

## 1. Introduction

Nonlinear evolution equations (NLEEs) represent a rich field for researchers due to their enormous applications in different engineering aspects. Fluid dynamics, shallow water waves, ocean waves, plasma explosions and wave propagation in optical fibers are like drops in the ocean of NLEEs applications. During last few decades, continuous efforts are exerted by different researchers to analyze and investigate different models of NLEEs either by analytic or numerical techniques. Similarity transformations [1-6], Lax-pair [1, 7-10], transformed rational function [11-13], $\acute{G}/G$ [13, 14] and sin-cosine methods [15-17] are some of famous methods in solving NLEEs. One of the most common evolution equations is Burgers' equation. Various numerical algorithms have been exploited to numerically solve NLEEs especially Burgers' equation to achieve minimized errors with respect to analytical solutions. Finite difference and other modifications [18-24], finite element and B-spline finite element [25-27], spectral least squares method [28-30], variational iteration method [31, 32], Adomian-Pade technique[33], homotopy analysis [34] and automatic differentiation method [35] are examples of such numerical techniques. Moreover, miscellaneous numerical techniques have been employed either for Burgers' equation or other engineering applications such as boundary element techniques for cavitation of hydrofoils [36, 37], finite element for channel flow with obstacles [38], modified diffusion coefficient technique for convection diffusion equation [39] and differential quadrature for functionally graded nanobeams [40, 41].

The proposed technique is based on discretizing the spatial direction using Chebyshev function while the differential quadrature method is employed to attain the spatial derivatives. Backward differentiation formula (BDF) has been utilized to obtain time derivative. To validate the proposed technique, it has been applied to one and dimensional Burgers' equation in addition to coupled model of the same equation. The present manuscript is arranged as follows. The proposed technique and solution procedures is presented in section 2. The error analysis of the considered scheme is being illustrated in section 3. The numerical attained results are shown in section 4. The paper is terminated by conclusions remarks in section 5.

## 2. The proposed technique

In this section, the proposed technique is illustrated to solve different cases of Burgers' equation. One, two-dimensional and coupled equations of Burgers' equation is being investigated throughout this paragraph.

### 2.1 One-dimensional Burgers' equation

Consider one-dimensional Burgers' equation in the following form:

$$\frac{\partial u}{\partial t} + u \frac{\partial u}{\partial x} = \nu \frac{\partial^2 u}{\partial x^2}, \quad P_1 \leq x \leq P_2, \ t \in [0, T], \tag{2.1}$$

with the initial condition

$$u(x, 0) = u_0(x), \tag{2.2}$$

and the Dirichlet boundary conditions

$$u(P_1, t) = u(P_2, t) = 0; \quad t \in [0, T], \tag{2.3}$$

The technique starts by discretization in time domain by using $p^{\text{th}}$ order backward differentiation formula (BDF) [42]. Equation (2.1) can be rewritten as:

$$\frac{\partial u}{\partial t} = -u \frac{\partial u}{\partial x} + \nu \frac{\partial^2 u}{\partial x^2}, \quad P_1 \leq x \leq P_2, \ t \in [0, T], \tag{2.4}$$

The time interval [0, T] is to be divided into N steps $0 = t_0 \leq t_1 \leq \cdots \leq t_N = T$, with constant time step $\Delta t = T/N$ and $t_n = n * \Delta t$ for $n = 1, 2, \ldots, N$. Then, second order backward differentiation formula (BDF-2) [23] along *t* domain is applied. Equation (2.4) becomes:

$$u^{n+1} = \frac{4}{3} u^n - \frac{1}{3} u^{n-1} + \frac{2}{3} \Delta t \left( \nu u_{xx}^{n+1} - u^{n+1} u_x^{n+1} \right) \tag{2.5}$$

Obviously, applying equation (2.5) at each time step requires obtaining the solutions in two previous time levels. Hence:

At *n=1*, the solution at time level $n - 1$ can be directly obtained using the initial condition while the solution at level $n$ can be attained using backward Euler Formula (BDF-1):

$$u^n = u^{n-1} + \Delta t \left( \nu u_{xx}^{n-1} - u^{n-1} u_x^{n-1} \right) \tag{2.6}$$

Now, for $n > 1$, the previous two-time levels are now ready to be used while the solution at the new time level, $n + 1$, can be attained using (2.5).

The nonlinear term appearing in (2.5) is required to be linearized to be ready for calculations. The linearization process is based on the approximation [43] $u^{n+1}u_x^{n+1} \approx w^{n+1}u_x^{n+1}$, where $w^{n+1}$ is computed by linear extrapolation using $u^n$ and $u^{n-1}$ as:

$$u^{n+1} \cong w^{n+1} = \left(1 + \left(\frac{K_{n+1}}{K_n}\right)\right)u^n - \left(\frac{K_{n+1}}{K_n}\right)u^{n-1} \tag{2.7}$$

where $K_{n+1} = t_{n+1} - t_n$ and $K_n = t_n - t_{n-1}$. Using (2.7), equation (2.5) is rewritten as:

$$u^{n+1} = \frac{4}{3}u^n - \frac{1}{3}u^{n-1} + \frac{2}{3}\Delta t \left(v\, u_{xx}^{n+1} - \left(\left(1 + \left(\frac{K_{n+1}}{K_n}\right)\right)u^n - \left(\frac{K_{n+1}}{K_n}\right)u^{n-1}\right)u_x^{n+1}\right) \tag{2.8}$$

Now, the discretization process in time domain has been accomplished. The next step in the proposed scheme is to discretize the spatial direction using Chebyshev function and use generalized differential quadrature method (GDQM) to construct the first and second derivatives [44-46]. Divide the interval $[P_1, P_2]$ into M nodes $P_1 = x_1 \leq x_2 \leq \cdots \leq x_M = P_2$ according to Chebyshev-Gauss-Lobatto distribution $x_i = x_1 + \frac{1}{2}(x_M - x_1)\left(1 - \cos\left(\frac{i-1}{M-1}\pi\right)\right), i = 1,2,\ldots M$. The derivatives are obtained using GDQM at a given node by approximating it using a weighted sum of the function values at all domain nodes according to the relation:

$$\left.\frac{d^m u}{dx^m}\right|_{x=x_i} = \sum_{j=1}^M a_{ij}^m u_j, \quad i = 1, 2, \ldots, M \tag{2.9}$$

where, $u_j = u(x_j)$ and $a_{ij}^m$ represent the corresponding weighting coefficients. The first derivative has a weighting coefficient in the form [44]:

$$a_{ij}^1 = \frac{1}{x_j - x_i}\left(\frac{Q_i}{Q_j}\right), i \neq j \quad and \quad a_{ii}^1 = -\sum_{j=1, i \neq 1}^M a_{ij}^1 \tag{2.10}$$

where

$$Q_i = \prod_{j=1, j \neq i}^M (x_i - x_j), \quad i, j = 1, 2, \ldots, M \tag{2.11}$$

To construct matrix form, the discrete values of dependent variable $u_i = u(x_i)$ at nodes are to be given as a vector $u = [u_1, u_2, \ldots u_M]^T$. Moreover, the first derivative vector is assumed to be $U^{(1)} = U$, then

$$U = Au \tag{2.12}$$

where, $A = [a_{ij}^1]$ represents the weighting coefficient matrix of first derivative. Matrix multiplication can be used to construct weighting coefficients for higher derivatives. For example, for the second derivative, $U^{(2)} = A\,U^{(1)} = AAu = A^{(2)}u$ so $A^{(2)} = A^2$ is weighting matrix for second derivative. Generally, for $m^{th}$ derivative, the coefficient matrix is given by $A^{(m)} = A^m$. Now, equation (2.8) is transformed to:

$$u^{n+1} = \frac{4}{3}u^n - \frac{1}{3}u^{n-1} + \frac{2}{3}\Delta t \left(v A^{(2)} u^{n+1} - \left(\left(1 + \left(\frac{K_{n+1}}{K_n}\right)\right)u^n - \left(\frac{K_{n+1}}{K_n}\right)u^{n-1}\right)A^{(1)} u^{n+1}\right) \quad (2.13)$$

In the proposed scheme, the time step is considered as a constant value. So, $K_{n+1} = K_n = \Delta t$ and equation (2.13) will be rearrange as:

$$\left[I - \frac{2}{3}\Delta t\, v\, A^{(2)} + \frac{2}{3}\Delta t (2u^n - u^{n-1})A^{(1)}\right] u^{n+1} = \frac{4}{3}u^n - \frac{1}{3}u^{n-1}. \quad (2.14)$$

Equation (2.14) can be written in a simple form:

$$Cu^{n+1} = F \quad (2.15)$$

where,

$$C = \left[I - \frac{2}{3}\Delta t\, v\, A^{(2)} + \frac{2}{3}\Delta t (2u^n - u^{n-1})A^{(1)}\right] \quad \text{and} \quad F = \frac{4}{3}u^n - \frac{1}{3}u^{n-1}.$$

Equation (2.15) can be solved efficiently by Thomas algorithm to get the solution vector $u = [u_1, u_2, \ldots u_M]$ at the step time (n+1).

*2.2 Two-dimensional Burgers' equation*

Consider two-dimensional Burgers' equation in the following form:

$$\frac{\partial u}{\partial t} + u\frac{\partial u}{\partial x} + u\frac{\partial u}{\partial y} = v\left(\frac{\partial^2 u}{\partial x^2} + \frac{\partial^2 u}{\partial y^2}\right), \quad z_1 \leq x \leq z_2,\ z_3 \leq y \leq z_4,\ t \in [0, T], \quad (2.16)$$

with the initial condition

$$u(x, y, 0) = u_0(x, y), \quad (2.17)$$

and the Dirichlet boundary conditions

$$u(z_1, y, t) = g_1(y, t), \quad u(z_2, y, t) = g_2(y, t), \quad u(x, z_3, t) = g_3(x, t),$$

$$u(x, z_4, t) = g_4(x, t); \quad t \in [0, T], \quad (2.18)$$

Equation (2.18) is discretized in time domain by the same proposed scheme described in section 2.1 with constant time step to be in the form:

$$u^{n+1} = \frac{4}{3}u^n - \frac{1}{3}u^{n-1} + \frac{2}{3}\Delta t(vu^{n+1}_{xx} + vu^{n+1}_{yy} - w^{n+1}u^{n+1}_x - w^{n+1}u^{n+1}_y) \quad (2.19)$$

where;

$$w^{n+1} = 2u^n - u^{n-1} \quad (2.20)$$

Equation (2.19) is discretized in the spatial directions and the derivatives are obtained as in section 2.1 to get the following equation:

$$u^{n+1} = \frac{4}{3}u^n - \frac{1}{3}u^{n-1} + \frac{2}{3}\Delta t \left(v\,[A^{(2)} + B^{(2)}]u^{n+1} - (2u^n - u^{n-1})[A^{(1)} + B^{(1)}]u^{n+1}\right) \quad (2.21)$$

where, $B^{(1)}$ and $B^{(2)}$, are the first and second derivatives respectively, which are obtained as in section 2.1 by applying the GDQM in y direction.

Equation (2.21) will be rearrange and written in a simple form:

$$Cu^{n+1} = F \quad (2.22)$$

where,

$$C = \left[I - \frac{2}{3}\Delta t\,v\,[A^{(2)} + B^{(2)}] + \frac{2}{3}\Delta t(2u^n - u^{n-1})[A^{(1)} + B^{(1)}]\right] \quad \text{and} \quad F = \frac{4}{3}u^n - \frac{1}{3}u^{n-1}.$$

Equation (2.22) can be solved efficiently by Thomas algorithm to get the solution vector $u = [u_1, u_2, \ldots u_M]$ at the step time (n+1).

*2.3 Two-dimensional Coupled Burgers' equations*:

Consider the following coupled equation

$$\frac{\partial u}{\partial t} + u\frac{\partial u}{\partial x} + v\frac{\partial u}{\partial y} - \frac{1}{Re}\left(\frac{\partial^2 u}{\partial x^2} + \frac{\partial^2 u}{\partial y^2}\right) = 0,$$

$$\quad (2.23)$$

$$\frac{\partial v}{\partial t} + u\frac{\partial v}{\partial x} + v\frac{\partial v}{\partial y} - \frac{1}{Re}\left(\frac{\partial^2 v}{\partial x^2} + \frac{\partial^2 v}{\partial y^2}\right) = 0,$$

subjected to the initial conditions:

$$u(x, y, 0) = f_1(x, y); (x, y) \in \Omega,$$

$$\quad (2.24)$$

$$v(x, y, 0) = f_2(x, y); (x, y) \in \Omega,$$

and boundary conditions:

$$u(x, y, t) = g_1(x, y); (x, y) \in \partial\Omega, t > 0,$$

$$\quad (2.25)$$

$$v(x, y, t) = g_2(x, y); (x, y) \in \partial\Omega, t > 0,$$

Equation (2.23) is discretized the time region by the same proposed scheme described in section 2.1 to be in the following form:

$$u^{n+1} = \frac{4}{3}u^n - \frac{1}{3}u^{n-1} + \frac{2}{3}\Delta t\left(\nu u_{xx}^{n+1} + \nu u_{yy}^{n+1} - w^{n+1}u_x^{n+1} - \eta^{n+1}u_y^{n+1}\right)$$

(2.26)

$$v^{n+1} = \frac{4}{3}v^n - \frac{1}{3}v^{n-1} + \frac{2}{3}\Delta t\left(\nu v_{xx}^{n+1} + \nu v_{yy}^{n+1} - w^{n+1}v_x^{n+1} - \eta^{n+1}v_y^{n+1}\right)$$

where;

$$w^{n+1} = 2u^n - u^{n-1} \quad , \quad \eta^{n+1} = 2v^n - v^{n-1} \tag{2.27}$$

As in the previous section 2.2, equation (2.26) will be written in a simple form:

$$Cu^{n+1} = F \ , \ Cv^{n+1} = G \tag{2.28}$$

where,

$$C = \left[I - \frac{2}{3}\Delta t \, \nu \, [A^{(2)} + B^{(2)}] + \frac{2}{3}\Delta t[(2u^n - u^{n-1})A^{(1)} + (2v^n - v^{n-1})B^{(1)}]\right],$$

$$F = \frac{4}{3}u^n - \frac{1}{3}u^{n-1}, \quad G = \frac{4}{3}v^n - \frac{1}{3}v^{n-1}.$$

Equation (2.28) can be solved efficiently by Thomas algorithm to get the solution vector $u = [u_1, u_2, \ldots u_M]$ and $v = [v, v_2, \ldots v_M]$ at the step time (n+1).

## 3. Stability Analysis

After discretization process using the proposed scheme, equation (2.1) are transformed to the following set of ordinary differential equations in time:

$$\frac{d\{U\}}{dt} = P\{U\} + \{E\} \tag{3.1}$$

while the coupled version of Burgers' equations (2.23) is transformed to:

$$\frac{d\{W\}}{dt} = R\{W\} + \{K\} \tag{3.2}$$

where;

R=$\begin{bmatrix} A & O \\ O & B \end{bmatrix}$,

and

1- $\{U\} = (u_2, u_3, \ldots, u_{N-1})$ is the unknown vector variable at the internal node of the grid.

2- $\{E\}$ is a vector containing non-homogenous part and boundary conditions.

3- $P = -\alpha U_{ij}A_1 + \nu A_2$

4- $O's$ are the null matrices.

5- $\{W\} = (U, V)^T$ are the unknown variables at the internal nodes of the grid:

$$\{U\} = (u_{22}, u_{23}, \ldots, u_{32}, u_{33}, \ldots, u_{42}, u_{43}, \ldots, u_{(N-1)(N-1)}),$$

$$\{V\} = (v_{22}, v_{23}, \ldots, v_{32}, v_{33}, \ldots, v_{42}, v_{43}, \ldots, v_{(M-1)(M-1)})$$

6- $\{K\}$ is the a vector containing non-homogenous part and boundary conditions.

7- $A = -U_{ij}A_1 - V_{ij}B_1 + vA_2 + vB_2$ and $B = -U_{ij}A'_1 - V_{ij}B'_1 + vA'_2 + vB'_2$.

$P_r$ are the matrices of the weighting coefficients $a_{ij}^{(r)}$ ($r = 1, 2$) respectively and given by:

$$P_r = \begin{bmatrix} a_{22}^{(r)} & a_{23}^{(r)} & \cdots & a_{2(N-1)}^{(r)} \\ a_{32}^{(r)} & a_{33}^{(r)} & & a_{3(N-1)}^{(r)} \\ \vdots & & \ddots & \vdots \\ a_{(N-1)2}^{(r)} & a_{(N-1)3}^{(r)} & \cdots & a_{(N-1)(N-1)}^{(r)} \end{bmatrix}_{(N-2) \times (N-2)},$$

while, $A_r$ and $B_r$ are square block diagonal matrices $(N-2) \times (M-2)$ of the weighting coefficients $a_{ij}^{(r)}$, $b_{ij}^{(r)}$ ($r = 1, 2$) respectively and given by:

$$A_r = \begin{bmatrix} a_{22}^{(r)}I & a_{23}^{(r)}I & \cdots & a_{2(N-1)}^{(r)}I \\ a_{32}^{(r)}I & a_{33}^{(r)}I & & a_{3(N-1)}^{(r)}I \\ \vdots & \vdots & \ddots & \vdots \\ a_{(N-1)2}^{(r)}I & a_{(N-1)3}^{(r)}I & \cdots & a_{(N-1)(N-1)}^{(r)}I \end{bmatrix},$$

$$B_r = \begin{bmatrix} M_r & 0 & \cdots & 0 \\ 0 & M_r & & 0 \\ \vdots & \vdots & \ddots & \vdots \\ 0 & 0 & \cdots & M_r \end{bmatrix}$$

where;

$$M_r = \begin{bmatrix} b_{22}^{(r)} & b_{23}^{(r)} & \cdots & b_{2(M-1)}^{(r)} \\ a_{32}^{(r)} & a_{33}^{(r)} & & b_{3(M-1)}^{(r)} \\ \vdots & \vdots & \ddots & \vdots \\ b_{(M-1)2}^{(r)} & b_{(M-1)3}^{(r)} & \cdots & b_{(M-1)(M-1)}^{(r)} \end{bmatrix}.$$

$I$ and $O$ are the matrices of order $(N-2)$ and $(M-2)$.
Also, $A'_r$ and $B'_r$ are square block diagonal matrices ( each of order $(N-2) \times (M-2)$) with the weighting coefficients $a_{ij}^{(r)}$, $b_{ij}^{(r)}$ ($r = 1, 2$) respectively.

    The stability of the proposed scheme depends on the stability of the system (3.1) and (3.2). The eigen values of the coefficient matrices P, A and B represent the base stone of stability analysis. The system (3.1) and (3.2) will be stable if the real part of each eigen value of P and R are either negative or zero which is illustrated by Figs. 1-3.

## 4. Results and discussion

$L_2$ and $L_\infty$ error norms are computed after each time step by using the following definitions:

$$L_2 := \|u_{exact} - u_{computed}\|_2 = \sqrt{\frac{\left(\sum_{j=1}^{n}\left|u_j^{exact} - u_j^{computed}\right|^2\right)}{N}}$$

(4.1)

$$L_\infty := \|u_{exact} - u_{computed}\|_\infty = max_j\left|u_j^{exact} - u_j^{computed}\right|$$

where $u_{exact}, u_{computed}$ are the exact and computed solutions at each point, respectively.

*4.1 Error analysis for 1-D Burgers' equation*

The analysis has been executed on two cases of relevant boundary and initial conditions.

*Case 1:*

Burgers' equation, (2.1)-(2.3), has been solve with initial and boundary conditions which are taken from the exact solution given by Wood [47]:

$$u(x,t) = f(x) = \frac{2\nu\pi\ e^{-\nu\pi^2 t}\ sin(\pi x)}{\sigma + e^{-\nu\pi^2 t}\ cos(\pi x)}, \quad \text{for } 0 \leq x \leq 1 \text{ and } t \geq 0 \quad \text{where, } \sigma > 1. \quad (4.2)$$

A comparison between numerical results of the proposed technique with Mittal and Jain work [48] and exact solution, equation (4.2), of case1 has been established with $\sigma=2$, N = 40, $\Delta t$ = 0.0001and
T= 0.001for Re=1, 10 at different nodes. The results are tabulated in Table 1. Another Comparison of $L_2$ and $L_\infty$ errors of case 1 with $\sigma=100$, T=1 and $\Delta t=0.01$ at different *Re* and N is illustrated in Table 2 for the present work compared to previous works [48-50]. Also, $L_2$ and $L_\infty$ errors of case 1 compared to Jiwari's work [51], with $\sigma=2$ and $\Delta t$ =0.001 for T= 0.1, 0.5 at different *Re*, are shown in Table 3. Figures 4 and 5 illustrate the numerical solution graphically at different values of $\sigma$ and *Re*, respectively.

*Case 2:*

In this case, Burgers' equations, (2.1)-(2.3), is solved with initial condition:

$$u(x, 0) = 4x(1 - x), 0 \leq x \leq 1 \quad (4.3)$$

and the boundary conditions:

$$u(0, t) = u(1, t) = 0, \ 0 \leq t \leq T \quad (4.4)$$

where exact solution is given by:

$$u(x,t) = 2\nu\pi\frac{\sum_{n=1}^{\infty} c_n \exp(-n^2\pi^2 \nu t)\ n\ sin(n\pi x)}{c_0 + \sum_{n=1}^{\infty} c_n \exp(-n^2\pi^2 \nu t)\ cos(n\pi x)} \quad (4.5)$$

and the Fourier coefficients $c_0$ and $c_n$ are, as in [22]:

$$c_0 = \int_0^1 \exp\left\{-\frac{1}{3v}[x^2(3-2x)]\right\} dx$$

$$c_n = 2\int_0^1 \exp\left\{-\frac{1}{3v}[x^2(3-2x)]\right\} \cos(n\pi x)\, dx$$

A comparison between the proposed numerical results with respect to other numerical Refs. [49, 51, 52] and the exact solution, equation (4.5), of case2 has been illustrated by Table 4 for Re =100 and Δt = 0.001 at different time and nodes. Table 5 indicates a comparison of $L_2$ and $L_\infty$ errors of case 2 for various time and viscosity at Δt =0.001, N=80. Figures, 6-8, illustrate graphical representation of the obtained numerical solutions at different values of $v$ and $Re$.

*4.2 Error analysis for 2-D Burgers' equation*

The solution of 2D Burgers' equation described by (2.16)-(2.18) is considered over a square domain $\Omega : [0,1] \times [0,1]$ with initial condition:

$$u(x,y,0) = \frac{1}{1+e^{\frac{Re(x+y)}{2}}}, \qquad (x,y) \in \Omega, \qquad (4.6)$$

and the boundary conditions:

$$u(0,y,t) = \frac{1}{1+e^{\frac{Re(y-t)}{2}}}, \qquad u(1,y,t) = \frac{1}{1+e^{\frac{Re(1+y-t)}{2}}}, \qquad y \in [0,1],\ t > 0,$$

$$u(x,0,t) = \frac{1}{1+e^{\frac{Re(x-t)}{2}}}, \qquad u(x,1,t) = \frac{1}{1+e^{\frac{Re(1+x-t)}{2}}}, \qquad x \in [0,1],\ t > 0, \quad (4.7)$$

for which the exact solution is given as:

$$u(x,y,t) = \frac{1}{1+e^{\frac{Re(x+y-t)}{2}}}, \quad (x,y) \in \Omega,\ t \geq 0. \qquad (4.8)$$

Table 6 illustrates absolute errors comparison between numerical solutions of Liu et al. [53] and the proposed scheme with Δt=0.001, Re =20 and grid size= 16×16 at some specific points for 2D Burgers' model. Also, the proposed scheme $L_2$ and $L_\infty$ error norms are compared to other authors' works [53-55] and the comparison is tabulated in Table 7. Moreover, Table 8 presents $L_2$ and $L_\infty$ error norms with Δt=0.0005 and grid size of 16×16 at different time and $Re$ values. Figure 9 depicts the numerical solution at different values of T, while Fig. 10 illustrates a comparison between numerical and exact solution of the equation.

*4.3 Error analysis for 2D coupled Burgers' equations*

The two-dimensional coupled Burgers' equations, (2.23) - (2.25), have been solve over a square domain $\Omega : [0,1] \times [0,1]$ with initial and boundary conditions considered according to exact solutions given by Li et al. work [56]:

$$u(x,y,t) = \frac{3}{4} - \frac{1}{4\left(1+e^{\frac{Re(4y-4x-t)}{32}}\right)} \ , \ v(x,y,t) = \frac{3}{4} + \frac{1}{4\left(1+e^{\frac{Re(4y-4x-t)}{32}}\right)} \ , (x,y) \in \Omega \ , t \geq 0. \quad (4.9)$$

Table 9 illustrates $L_2$ and $L_\infty$ error norms for u component at Re =100 and different time. A comparison with respect to the exact solution, equation (4.9), and previous works of Shukla et al [57] and Shi et al [58] of u is illustrated in Table 10 and v in Table 11 at Re =100 and 20×20 grid size. A graphical illustration of both numerical and exact solution is shown in Fig. 11.

## 5. Conclusions

The proposed numerical scheme is based on second order backward differentiation formula (BDF2) in time space and generalized differential quadrature method (GDQM) in spatial space. The scheme has been employed to solve Burgers' equation in three cases, one-dimensional, two-dimensional and coupled two-dimensional models. The numerical results showed the effectiveness of the proposed scheme compared with the numerical schemes in the literature. The obtained errors in all cases are less than those obtained by other schemes at different values of *Re* and $v$. Tables2 and 3 show that the proposed scheme produces error norms, $L_2$ and $L_\infty$, of significant decrease percentage compared to previous works of Mittal and Jain [**48**], Rahman et al. [**49**], Tamsir et al. [**50**] and Jiwari [**51**]. Moreover, in case of (2+1)-dimensional model, Tables 6 and 7 illustrate another superiority of the proposed scheme compared to other works of Liu et al. [**53**], Mittal and Tripathi [**54**] and Arora and Joshi [**55**]. Besides, in case of (2+1)-dimensional coupled model, the proposed scheme emphasizes its efficiency in reducing error norms compared to Li et al work [**56**] at different times. The stability analysis has been also executed emphasizing the stability of the proposed technique on the same manner described by Tamsir et al. [**50, 59**].

**References**.

[1]     S. M. Mabrouk and A. S. Rashed, "*Analysis of (3 + 1)-dimensional Boiti – Leon –Manna–Pempinelli equation via Lax pair investigation and group transformation method*", Computers & Mathematics with Applications **74** (2017), 2546-2556.

[2]     A. S. Rashed, "*Analysis of (3+1)-dimensional unsteady gas flow using optimal system of Lie symmetries*", Mathematics and Computers in Simulation **156** (2019), 327-346.

[3]     A. S. Rashed and M. M. Kassem, "*Group analysis for natural convection from a vertical plate*", Journal of Computational and Applied Mathematics **222** (2008), 392-403.

[4]     M. M. Kassem and A. S. Rashed, "*Group solution of a time dependent chemical convective process*", Applied Mathematics and Computation **215** (2009), 1671-1684.

[5]     A. S. Rashed and M. M. Kassem, "*Hidden symmetries and exact solutions of integro-differential Jaulent-Miodek evolution equation*", AMC Applied Mathematics and Computation **247** (2014), 1141-1155.

[6]     M. M. Kassem and A. S. Rashed, "*N-solitons and cuspon waves solutions of (2 + 1)-dimensional Broer–Kaup–Kupershmidt equations via hidden symmetries of Lie optimal system*", Chinese Journal of Physics **57** (2019), 90-104.

[7]     R. Saleh, M. Kassem and S. Mabrouk, "*Exact solutions of Calgero-Bogoyavlenskii-Schiff equation using the singular manifold method after Lie reductions*", MMA Mathematical Methods in the Applied Sciences **40** (2017), 5851-5862.

[8]     S. Mabrouk, M. Kassem and M. Abd-el-Malek, "*Group similarity solutions of the lax pair for a generalized HirotaSatsuma equation*", Applied Mathematics and Computation Applied Mathematics and Computation **219** (2013), 7882-7890.

[9]     P. G. Estévez, M. L. Gandarias and J. Prada, "*Symmetry reductions of a 2 + 1 Lax pair*", PLA Physics Letters A **343** (2005), 40-47.

[10]    M. Legare, "*Symmetry Reductions of the Lax Pair of the Four-Dimensional Euclidean Self-Dual Yang-Mills Equations*", JNMP Journal of Non-linear Mathematical Physics **3** (2007), 266.

[11]    S. M. Mabrouk and A. S. Rashed, "*N-Solitons, Kink and Periodic Wave Solutions for (3+1)-Dimensional Hirota Bilinear Equation Using Three Distinct Techniques*", Chinese Journal of Physics Chinese Journal of Physics (2019), 793.

[12]    W. X. Ma and J. H. Lee, "*A transformed rational function method and exact solutions to the 3 + 1 dimensional Jimbo-Miwa equation*", Chaos Solitons Fractals Chaos, Solitons and Fractals **42** (2009), 1356-1363.

[13]    S. Zhang and H.-Q. Zhang, "*A transformed rational function method for (3+1)-dimensional potential Yu-Toda-Sasa-Fukuyama equation*", Pramana - J Phys Pramana : Journal of Physics **76** (2011), 561-571.

[14]    H.-L. Lü, X.-Q. Liu and L. Niu, "*A generalized (G'/G) expansion method and its applications to nonlinear evolution equations*", Applied Mathematics and Computation **215** (2010), 3811–3816.

[15]    M. Yaghobi Moghaddam, A. Asgari and H. Yazdani, "*Exact travelling wave solutions for the generalized nonlinear Schrödinger (GNLS) equation with a source by Extended tanh–coth, sine–cosine and Exp-Function methods*", Applied Mathematics and Computation **210** (2009), 422-435.

[16]    G. Betchewe, B. B. Thomas, K. K. Victor and K. T. Crepin, "*Explicit series solutions to nonlinear evolution equations: The sine–cosine method*", Applied Mathematics and Computation **215** (2010), 4239-4247.

[17]    E. M. Zayed and M. A. Abdelaziz, "*Exact solutions for the nonlinear Schrödinger equation with variable coefficients using the generalized extended tanh-function, the sine–cosine and the exp-function methods*", Applied Mathematics and Computation **218** (2011), 2259-2268.

[18]    M. Ciment, S. H. Leventhal and B. C. Weinberg, "*The operator compact implicit method for parabolic equations*", Journal of Computational Physics Journal of Computational Physics **28** (1978), 135-166.

[19]    L. Iskandar and A. Mohsen, "*Some numerical experiments on the splitting of Burgers' equation*", Numerical Methods for Partial Differential Equations **8** (1992), 267-276.

[20]    I. A. Hassanien, A. A. Salama and H. A. Hosham, "*Fourth-order finite difference method for solving Burgers equation*", Applied Mathematics and Computation **170** (2005), 781-800.

[21]    T. Zhanlav, O. Chuluunbaatar and V. Ulziibayar, "*Higher-order accurate numerical solution of unsteady Burgers equation*", AMC Applied Mathematics and Computation **250** (2015), 701-707.

[22]    V. Mukundan and A. Awasthi, "*Efficient numerical techniques for Burgers' equation*", Applied Mathematics and Computation Applied Mathematics and Computation **262** (2015), 282-297.

[23]    S. Kutluay, A. R. Bahadir and A. Özdes, "*Numerical solution of one-dimensional Burgers equation: explicit and exact-explicit finite difference methods*", Journal of Computational and Applied Mathematics **103** (1999), 251-261.

[24]    N. A. Mohamed, "*Fully Implicit Scheme for Solving Burgers' Equation based on Finite Difference Method*", The Egyptian International Journal of Engineering Sciences and Technology **26** (2018), 38-44.

[25]    T. Ozis, E. N. Aksan and A. Ozdes, "*A finite element approach for solution of Burgers `equation*", Applied Mathematics and Computation **139** (2003), 417-428.

[26]    A. Dogan, "*A Galerkin finite element approach to Burgers' equation*", Applied mathematics and computation. **157** (2005), 331.

[27] E. N. Aksan, "*Quadratic B-spline finite element method for numerical solution of the Burgers equation*", Applied Mathematics and Computation **174** (2006), 884-896.
[28] B. De Maerschalck and M. Gerritsma, "*The Use of Chebyshev Polynomials in the Space-Time Least-Squares Spectral Element Method*", Numerical Algorithms **38** (2005), 173-196.
[29] B. D. Maerschalck and M. I. Gerritsma, "*Least-squares spectral element method for non-linear hyperbolic differential equations*", Journal of Computational and Applied Mathematics Journal of Computational and Applied Mathematics **215** (2008), 357-367.
[30] W. Heinrichs, "*An adaptive spectral least-squares scheme for the Burgers equation*", Numer Algor Numerical Algorithms **44** (2007), 1-10.
[31] M. A. Abdou and A. A. Soliman, "*Variational iteration method for solving Burger's and coupled Burger's equations*", CAM Journal of Computational and Applied Mathematics **181** (2005), 245-251.
[32] J. Biazar and H. Aminikhah, "*Exact and numerical solutions for non-linear Burger's equation by VIM*", Math. Comput. Model. Mathematical and Computer Modelling **49** (2009), 1394-1400.
[33] M. Dehghan, A. Hamidi and M. Shakourifar, "*The solution of coupled Burgers equations using AdomianPade technique*", Applied Mathematics and Computation Applied Mathematics and Computation **189** (2007), 1034-1047.
[34] M. M. Rashidi, G. Domairry and S. Dinarvand, "*Approximate solutions for the Burger and regularized long wave equations by means of the homotopy analysis method*", Communications in Nonlinear Science and Numerical Simulation **14** (2009), 708-717.
[35] A. Asaithambi, "*Numerical solution of the Burgers equation by automatic differentiation*", AMC Applied Mathematics and Computation **216** (2010), 2700-2708.
[36] N. A. Mohmaed, A. Abdel Gawad and A. S. G, "A Boundary Element Approach with New Algorithm to Predict Partial Cavitation around 2D Hydrofoils," *ICFDP8: Eighth International Congress of Fluid Dynamics & Propulsion*, Sharm El-Shiekh, Egypt, 2006, pp. 14-17.
[37] N. A. Mohamed, A. A. Gawad and S. G. Ahmed, "*Computation of Partial Cavitation Characteristics over Two-Dimensional Symmetric Hydrofoils Using A Newly Proposed Boundary Element Algorithm*", International eJournal of Engineering Mathematics: Theory and Application **1** (2007).
[38] A. Abdel Gawad, N. A. Mohamed, S. A. Mohamed and M. S. Matbuly, "*Investigation of the Channel Flow with Internal Obstacles Using Large Eddy Simulation and Finite-Element Technique*", Applied and computational mathematics **2** (2013), 1-13.
[39] S. A. Mohamed, N. A. Mohamed, A. F. Abdel Gawad and M. S. Matbuly, "*A modified diffusion coefficient technique for the convection diffusion equation*", Applied Mathematics and Computation Applied Mathematics and Computation **219** (2013), 9317-9330.
[40] R. A. Shanab, M. Attia, s. Mohamed and N. Alaa, "*Analytical Solution for Bending of Functionally Graded Timoshenko Nanobeam Incorporating Surface Energy and Microstructure Effects*", East African Scholars Journal of Engineering and Computer Sciences **2** (2019), 91-96.
[41] M. A. Attia, R. A. Shanab, S. A. Mohamed and N. A. Mohamed, "*Surface Energy Effects on the Nonlinear Free Vibration of Functionally Graded Timoshenko Nanobeams Based on Modified Couple Stress Theory*", International Journal of Structural Stability and Dynamics (2019).
[42] K. E. Atkinson, W. Han and D. Stewart, "*Numerical solutions of ordinary differential equations*", J. Wiley, New York, 2009.
[43] D. A. Kay, P. M. Gresho, D. F. Griffiths and D. J. Silvester, "*Adaptive Time-Stepping for Incompressible Flow Part II: Navier-Stokes Equations*", SIAM journal on scientific computing : a publication of the Society for Industrial and Applied Mathematics. **32** (2010), 111.
[44] C. Shu, "*Differential quadrature and its application in engineering*", Springer-Verlag, New York, N.Y., 2000.

[45]   M. S. Matbuly, O. Ragb and M. Nassar, "*Natural frequencies of a functionally graded cracked beam using the differential quadrature method*", AMC Applied Mathematics and Computation **215** (2009), 2307-2316.

[46]   M. A. Attia and S. A. Mohamed, "*Nonlinear modeling and analysis of electrically actuated viscoelastic microbeams based on the modified couple stress theory*", APM Applied Mathematical Modelling **41** (2017), 195-222.

[47]   W. L. Wood, "*An exact solution for Burgers equation*", Communications in Numerical Methods in Engineering **22** (2006), 797-798.

[48]   R. C. Mittal and R. K. Jain, "*Numerical solutions of nonlinear Burgers equation with modified cubic B-splines collocation method*", AMC Applied Mathematics and Computation **218** (2012), 7839-7855.

[49]   K. Rahman, N. Helil, R. Yimin, A. International Conference on Computer and M. System, "*Some new semi-implicit finite difference schemes for numerical solution of Burgers equation*", **14** (2010), V14-451-V414-455.

[50]   M. Tamsir, V. K. Srivastava and R. Jiwari, "*An algorithm based on exponential modified cubic B-spline differential quadrature method for nonlinear Burgers equation*", AMC Applied Mathematics and Computation **290** (2016), 111-124.

[51]   R. Jiwari, "*A hybrid numerical scheme for the numerical solution of the Burgers equation*", Computer Physics Communications Computer Physics Communications **188** (2015), 59-67.

[52]   H. Nojavan, S. Abbasbandy and M. Mohammadi, "*Local variably scaled Newton basis functions collocation method for solving Burgers equation*", Appl. Math. Comput. Applied Mathematics and Computation **330** (2018), 23-41.

[53]   X. Liu, J. Wang and Y. Zhou, "*A space-time fully decoupled wavelet Galerkin method for solving two-dimensional Burgers equations*", CAMWA Computers and Mathematics with Applications **72** (2016), 2908-2919.

[54]   R. C. Mittal and A. Tripathi, "*Numerical solutions of two-dimensional Burgers' equations using modified Bi-cubic B-spline finite elements*", Eng. Comput. (Swansea Wales) Engineering Computations (Swansea, Wales) **32** (2015), 1275-1306.

[55]   G. Arora and V. Joshi, "*A computational approach using modified trigonometric cubic B-spline for numerical solution of Burgers equation in one and two dimensions*", Alexandria Engineering Journal Alexandria Engineering Journal **57** (2018), 1087-1098.

[56]   Q. Li, Z. Chai and B. Shi, "*Lattice Boltzmann models for two-dimensional coupled Burgers equations*", Computers & mathematics with applications **75** (2018), 864-875.

[57]   H. S. Shukla, M. Tamsir, V. K. Srivastava and J. Kumar, "*Numerical solution of two dimensional coupled viscous Burger equation using modified cubic B-spline differential quadrature method*", AIP Advances **4** (2014), 117134.

[58]   F. Shi, H. Zheng, Y. Cao, J. Li and R. Zhao, "*A fast numerical method for solving coupled Burgers' equations Fast Numerical Method*", Numerical Methods for Partial Differential Equations **33** (2017), 1823-1838.

[59]   M. Tamsir, N. Dhiman and V. K. Srivastava, "*Extended modified cubic B-spline algorithm for nonlinear Burgers' equation*", Beni-Suef University journal of basic and applied sciences **5** (2016), 244-254.

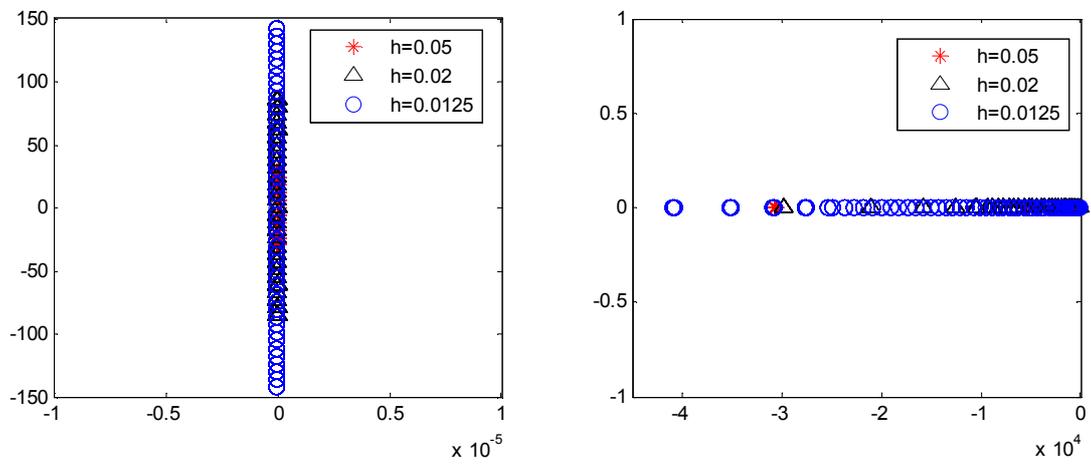

*Fig.1: Eigen values of $P_1$ (left) and $P_2$ (right) for different grid sizes.*

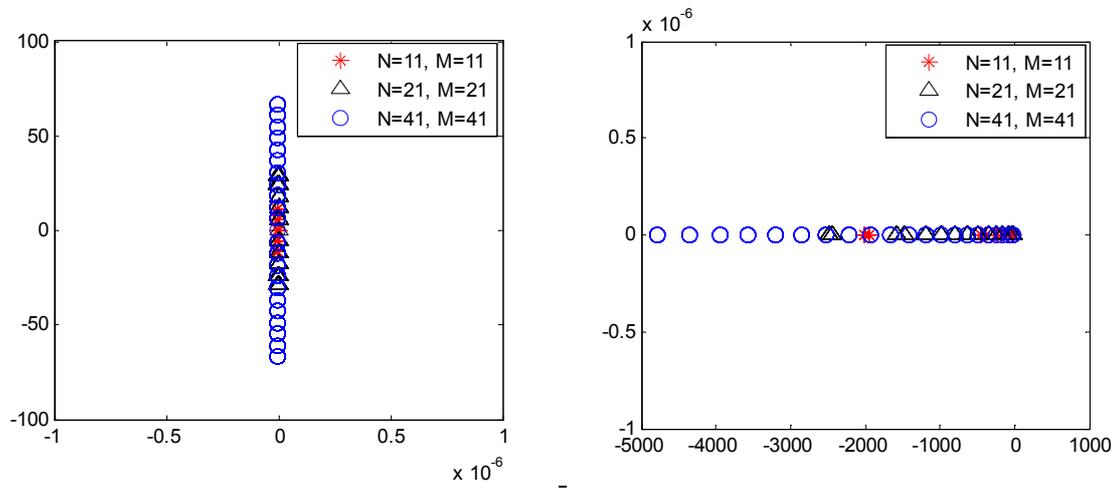

*Fig.2: Eigen values of $A_1$ (left) and $A_2$ (right) for different grid sizes.*

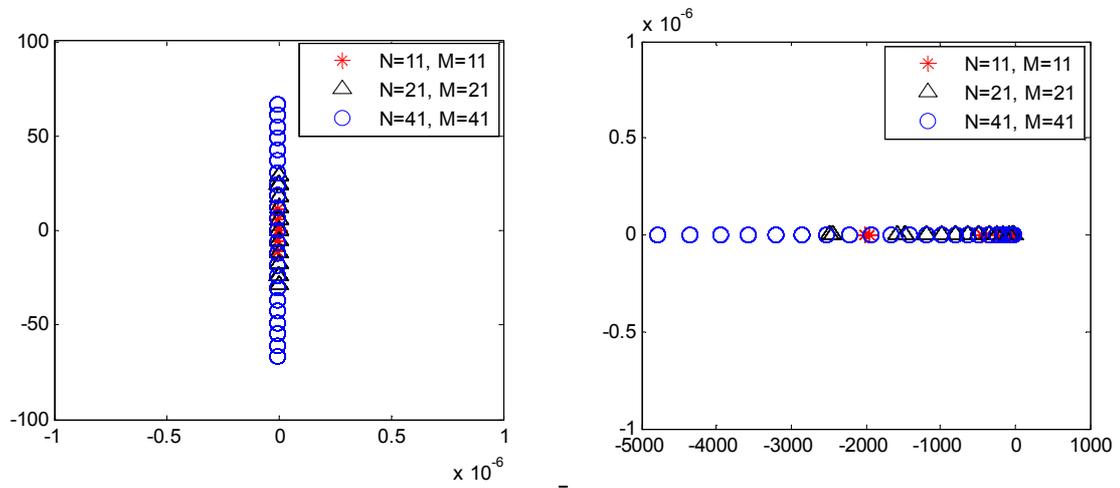

*Fig.3: Eigen values of $B_1$ (left) and $B_2$ (right) for different grid sizes.*

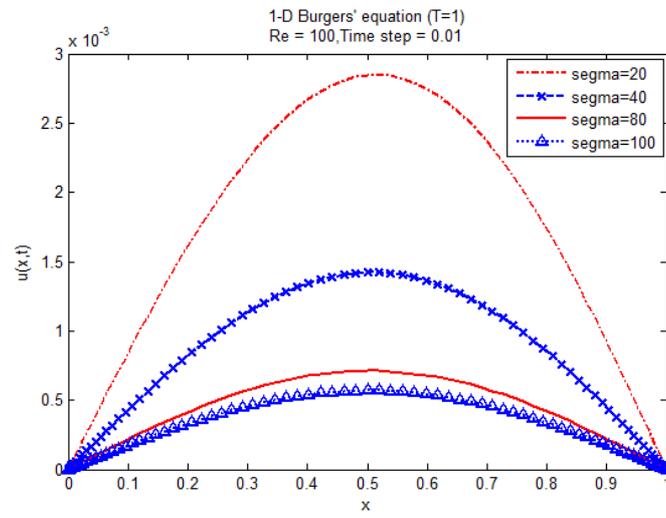

*Fig.4: Numerical solution at Re=100, Δt=0.01, T=1 and N=20 for different values of σ.*

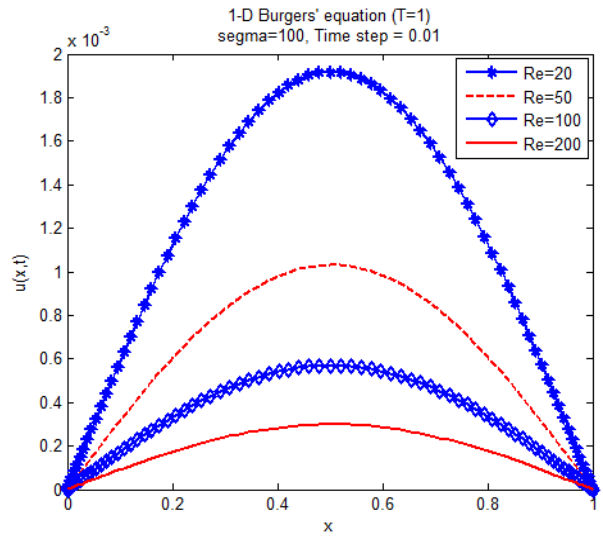

*Fig.5: Numerical solution at σ =100, Δt=0.01, T=1 and N=20 for different values of Re.*

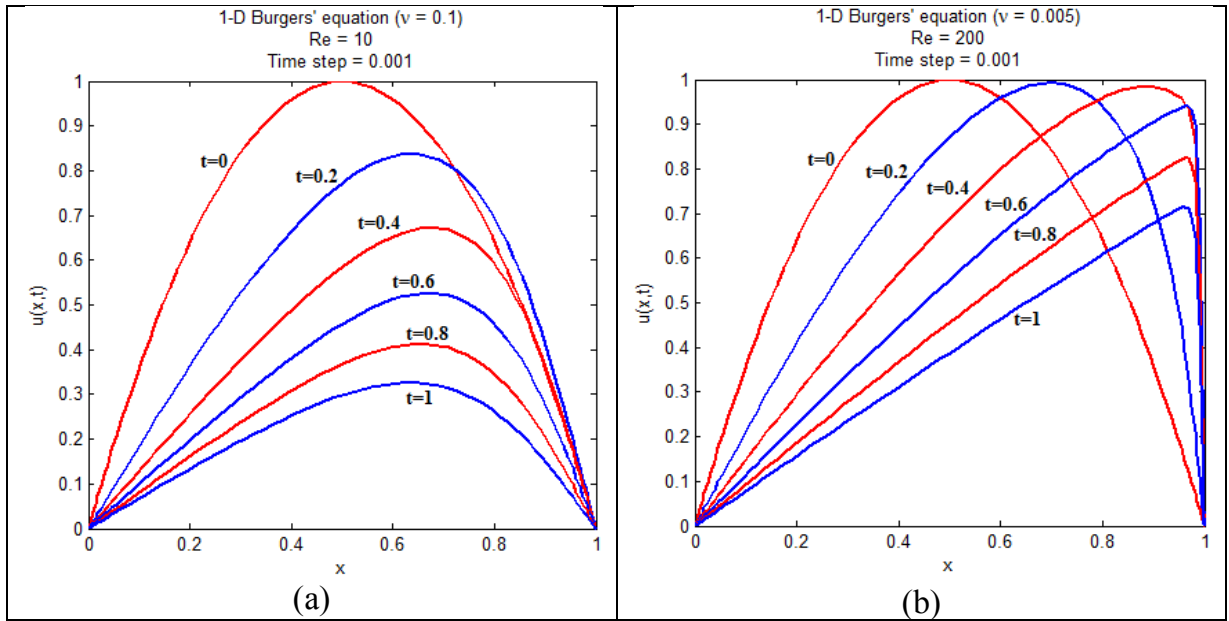

*Fig.6: Numerical solution of case 2 at different time t with Δt =0.001, N=80
for (a) Re=10 and (b) Re=200.*

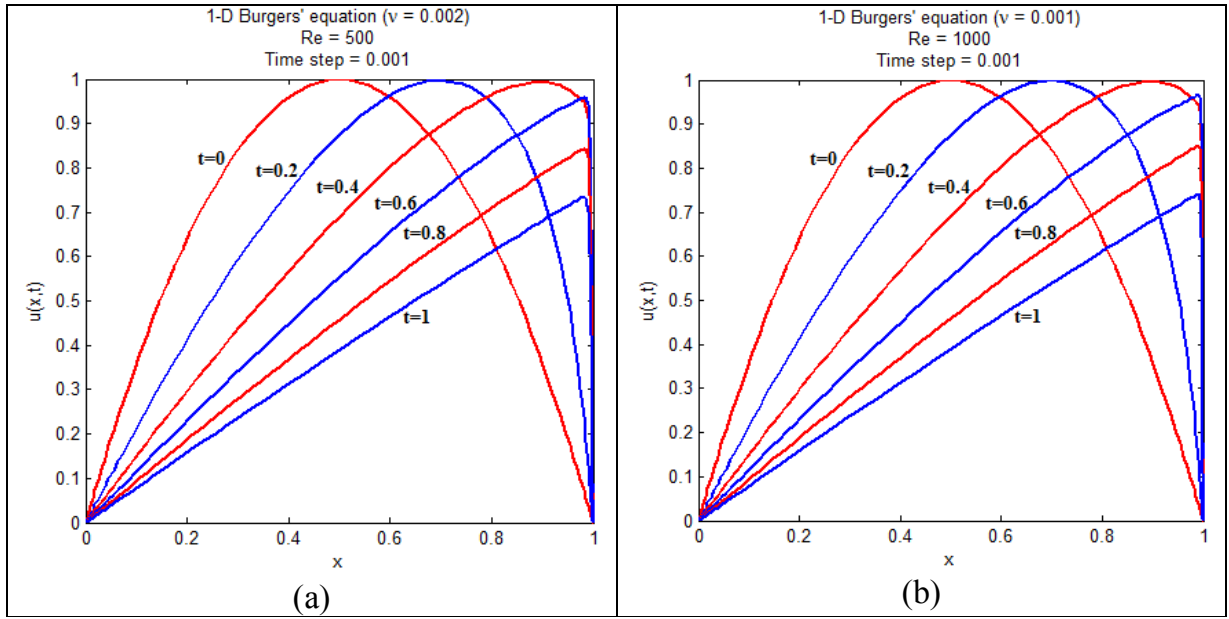

*Fig.7: Numerical solution of case 2 at different time t with Δt =0.001 for
(a) Re=500, N=80 (b) Re=1000, N=80.*

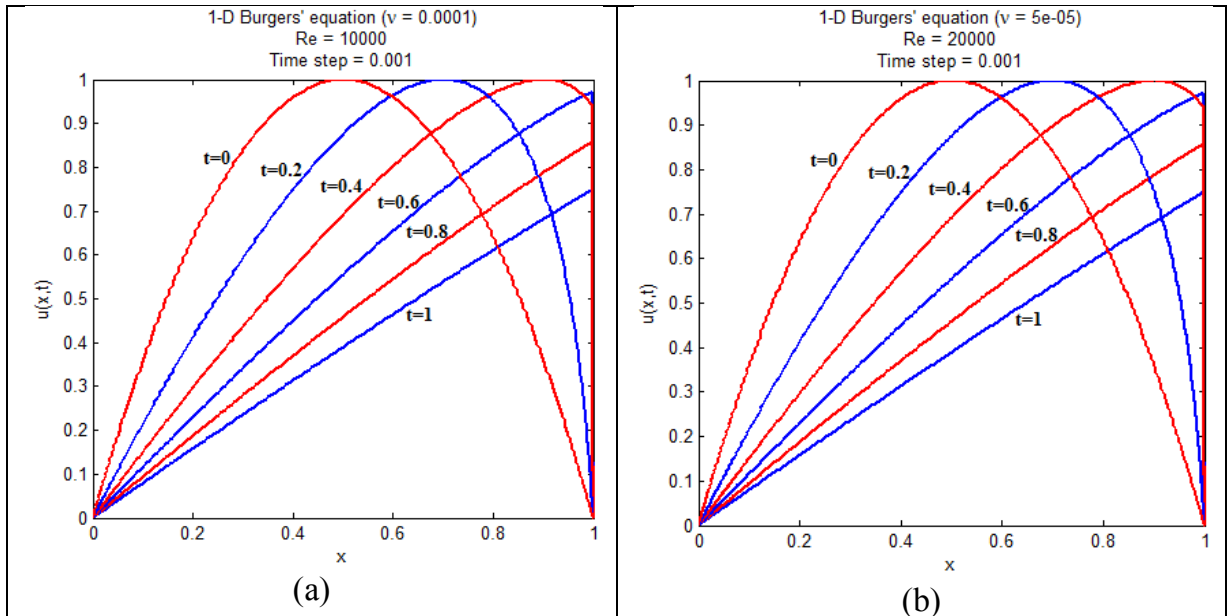

*Fig.8: Numerical solution of case 2 at different time t with Δt =0.001 for
(a) Re=10000, N=300  (b) Re=20000, N=400.*

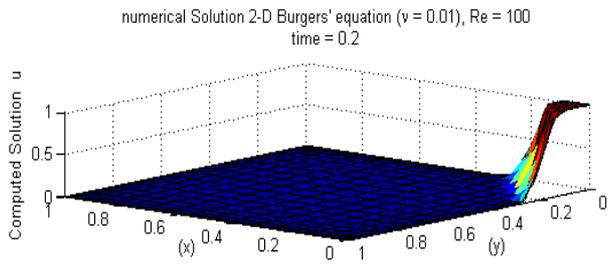 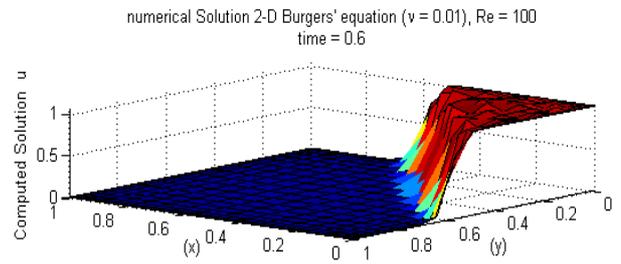
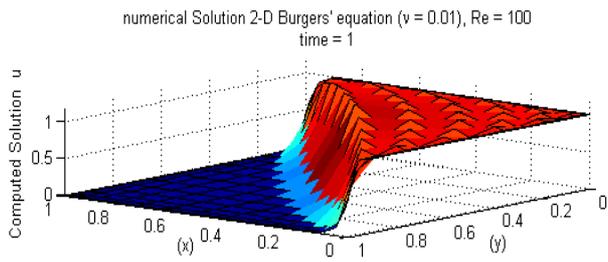 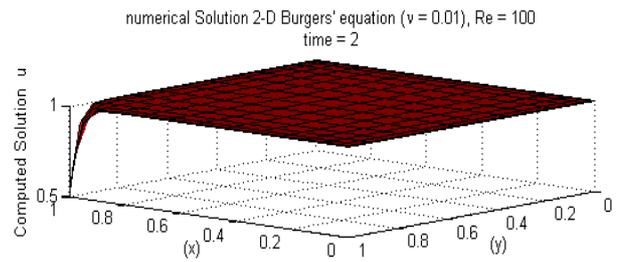

*Fig.9: Numerical solution of 2D Burgers' model at different time T = 0.2, 0.6, 1.0, 2.0 with∆t =0.001 for Re =100, grid size =16×16.*

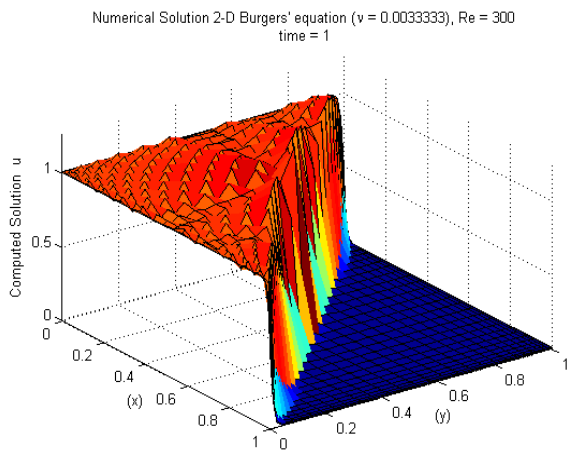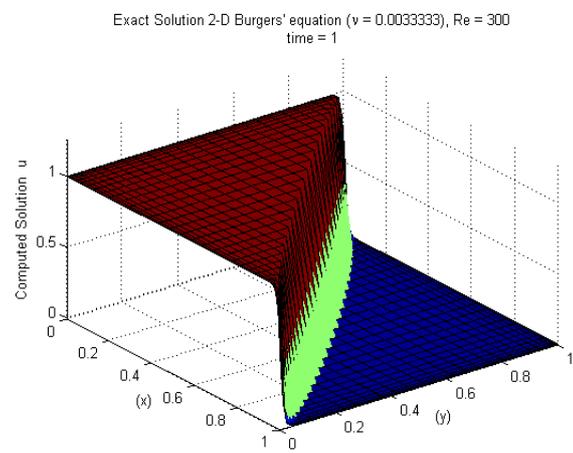

*Fig.10: Numerical and exact solutions of 2D Burgers' model at time T = 1 with Δt=0.0001 for Re =300, grid size =32×32.*

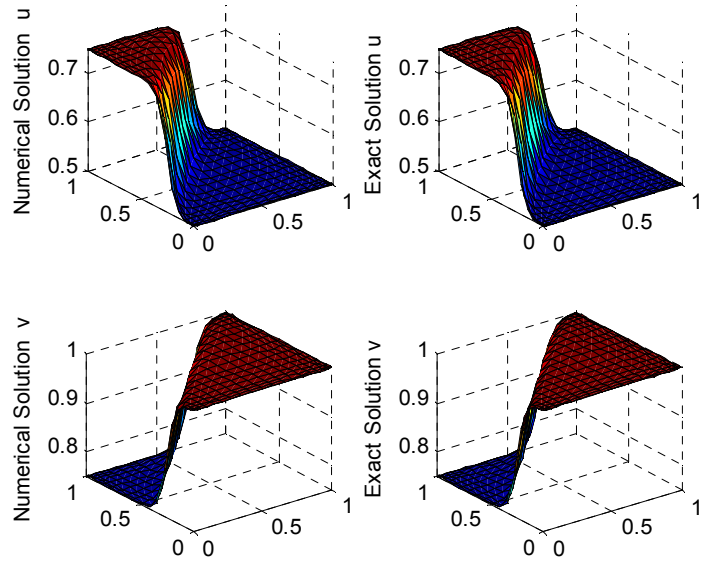

*Fig.11: Numerical and Exact solutions at T=1.0 with Re =200, Δt= 0.001 and grid size 20×20 for coupled Burgers' equation.*

*Table 1: comparison between the proposed numerical results with respect to other numerical Refs. and the exact solution case1 with σ=2, N = 40, Δt = 0.0001and T= 0.001 for Re=1, 10 at different nodes.*

| $x$ | $Re = 1$ | | | $Re = 10$ | | |
|---|---|---|---|---|---|---|
| | Mittal and Jain[48] | Proposed scheme | Exact solution | Mittal and Jain[48] | Proposed scheme | Exact solution |
| 0.1 | 0.653547 | 0.653544 | 0.653545 | 0.065750 | 0.065749761 | 0.065749761 |
| 0.2 | 1.305540 | 1.305533 | 1.305534 | 0.131383 | 0.131382943 | 0.131382943 |
| 0.3 | 1.949376 | 1.949363 | 1.949364 | 0.196281 | 0.196280911 | 0.196280911 |
| 0.4 | 2.565949 | 2.565927 | 2.565927 | 0.258576 | 0.258575994 | 0.258575995 |
| 0.5 | 3.110778 | 3.110739 | 3.110739 | 0.313850 | 0.313849356 | 0.313849356 |
| 0.6 | 3.492910 | 3.492873 | 3.492871 | 0.352972 | 0.352972354 | 0.352972351 |
| 0.7 | 3.549585 | 3.549602 | 3.549594 | 0.359443 | 0.359442750 | 0.359442742 |
| 0.8 | 3.049957 | 3.050145 | 3.050130 | 0.309579 | 0.309579979 | 0.309579963 |
| 0.9 | 1.816379 | 1.816672 | 1.816658 | 0.184751 | 0.184753526 | 0.184753511 |
| $L_2$ | 2.85 E-04 | 6.66 E-06 | -- | 3.08E-06 | 7.13 E-09 | -- |
| $L_\infty$ | 1.07 E-04 | 1.60 E-05 | -- | 1.15E-06 | 1.73 E-08 | -- |

*Table 2: Comparison of $L_2$ and $L_\infty$ errors of example 1 with $\sigma=100$, $T=1$ and $\Delta t=0.01$ at different Re and N.*

| N  | Re  | Rahman et al. [49] | | Mittal and Jain [48] | | Proposed scheme | |
|----|-----|------------|------------|------------|------------|------------|------------|
|    |     | $L_2$ | $L_\infty$ | $L_2$ | $L_\infty$ | $L_2$ | $L_\infty$ |
| 10 | 100 | 3.4545E-07 | 4.8808E-07 | 3.2840E-07 | 4.6280E-07 | 2.3494E-10 | 3.9698E-10 |
| 20 |     | 1.0124E-07 | 1.4305E-07 | 8.1921E-08 | 1.1640E-07 | 2.3486E-10 | 3.9784E-10 |
| 40 |     | 4.0028E-08 | 5.6677E-08 | 2.0470E-08 | 2.9068E-08 | 2.3486E-10 | 3.9784E-10 |
| 80 |     | 2.4713E-08 | 3.4992E-08 | 5.1194E-09 | 7.2706E-09 | 2.3486E-10 | 3.9856E-10 |

| N  | Re  | Mittal and Jain [48] | | Tamsir et al. [50] | | Proposed scheme | |
|----|-----|------------|------------|------------|------------|------------|------------|
|    |     | $L_2$ | $L_\infty$ | $L_2$ | $L_\infty$ | $L_2$ | $L_\infty$ |
| 10 | 200 | 8.631E-08  | 1.215E-07  | 6.330E-08  | 1.467E-07  | 3.161E-11  | 5.303E-11  |
| 20 |     | 2.153E-08  | 3.062E-08  | 1.014E-08  | 3.029E-08  | 3.155E-11  | 5.342E-11  |
| 40 |     | 5.378E-09  | 7.644E-09  | 1.207E-09  | 3.956E-09  | 3.155E-11  | 5.347E-11  |
| 80 |     | 1.345E-09  | 1.917E-09  | 1.322E-10  | 8.861E-11  | 3.155E-11  | 5.355E-11  |

*Table 3: Comparison of $L_2$ and $L_\infty$ errors of example 1 with $\sigma=2$ and $\Delta t=0.001$ for $T=0.1, 0.5$ at different Re.*

| T | Re | Jiwari [51] $\Delta t=0.001$, N=32 | | Proposed scheme $\Delta t=0.001$, N=20 | | Proposed scheme $\Delta t=0.001$, N=32 | |
|---|---|---|---|---|---|---|---|
| | | $L_2$ | $L_\infty$ | $L_2$ | $L_\infty$ | $L_2$ | $L_\infty$ |
| 0.1 | 10 | 3.1398E-07 | 3.5581E-06 | 3.6018E-07 | 8.1267E-07 | 3.6020E-07 | 8.1808E-07 |
| | 100 | 2.0149E-09 | 3.3607E-08 | 6.6321E-10 | 1.5636E-09 | 6.6393E-10 | 1.5970E-09 |
| | 10000 | 1.1730E-12 | 1.1040E-10 | 3.5865E-14 | 1.4146E-13 | 7.1838E-16 | 1.7435E-15 |
| | 100000 | 3.0739E-14 | 3.0451E-12 | 4.1011E-16 | 1.6284E-15 | 8.1422E-19 | 2.0939E-18 |
| 0.5 | 10 | 7.5254E-07 | 7.3454E-06 | 6.2788E-08 | 1.2138E-07 | 6.2787E-08 | 1.2158E-07 |
| | 100 | 7.2279E-09 | 1.1066E-07 | 4.9320E-10 | 1.1375E-09 | 4.9473E-10 | 1.1514E-09 |
| | 10000 | 2.1651E-12 | 1.4078E-10 | 1.0714E-13 | 4.0785E-13 | 7.0998E-16 | 1.7148E-15 |
| | 100000 | 8.8788E-14 | 8.7450E-12 | 1.9304E-15 | 7.6459E-15 | 3.3693E-18 | 1.1923E-17 |

*Table 4: comparison between the proposed numerical results with respect to other numerical Refs. and the exact solution of case 2 for Re =100 and Δt = 0.001 at different time and nodes.*

| x | t | Mittal [48] $\Delta t$=0.001 N=40 | Jiwari [51] $\Delta t$=0.001 | Nojavan [52] N=180 | Proposed scheme $\Delta t$=0.001 N=40 | Proposed scheme $\Delta t$=0.001 N=80 | Exact solution | Absolute error Nojavan [52] | Absolute error Proposed scheme |
|---|---|---|---|---|---|---|---|---|---|
| 0.25 | 0.4 | 0.36225 | 0.36225 | 0.36226 | 0.36226 | 0.36226 | 0.36226 | 6.1 E−5 | 9.9 E−8 |
|  | 0.6 | 0.28202 | 0.28204 | 0.28204 | 0.28204 | 0.28204 | 0.28204 | 5.1 E−6 | 7.0 E−8 |
|  | 0.8 | 0.23044 | 0.23045 | 0.23045 | 0.23045 | 0.23045 | 0.23045 | 7.6 E−6 | 4.6 E−8 |
|  | 1.0 | 0.19468 | 0.19469 | 0.19469 | 0.19469 | 0.19469 | 0.19469 | 1.1 E−6 | 3.2 E−8 |
|  | 3.0 | 0.07613 | 0.07613 | 0.07613 | 0.07613 | 0.07613 | 0.07613 | 4.1 E−6 | 4.7 E−9 |
| 0.5 | 0.4 | 0.68368 | 0.68364 | 0.68368 | 0.68369 | 0.68369 | 0.68368 | 2.6 E−6 | 1.2 E−7 |
|  | 0.6 | 0.54832 | 0.54831 | 0.54832 | 0.54832 | 0.54832 | 0.54832 | 1.4 E−5 | 1.7 E−7 |
|  | 0.8 | 0.45371 | 0.45371 | 0.45370 | 0.45372 | 0.45371 | 0.45371 | 3.1 E−6 | 1.4 E−7 |
|  | 1.0 | 0.38567 | 0.38568 | 0.38568 | 0.38568 | 0.38568 | 0.38568 | 7.0 E−6 | 9.7 E−8 |
|  | 3.0 | 0.15218 | 0.15219 | 0.15218 | 0.15218 | 0.15218 | 0.15218 | 6.6 E−6 | 1.1 E−8 |
| 0.75 | 0.4 | 0.92052 | 0.92044 | 0.92051 | 0.92050 | 0.92050 | 0.92050 | 1.9 E−5 | 2.5 E−7 |
|  | 0.6 | 0.78300 | 0.78297 | 0.78302 | 0.78299 | 0.78299 | 0.78299 | 3.8 E−5 | 4.3 E−7 |
|  | 0.8 | 0.66272 | 0.66272 | 0.66272 | 0.66272 | 0.66272 | 0.66272 | 8.1 E−6 | 3.2 E−7 |
|  | 1.0 | 0.56932 | 0.56932 | 0.56932 | 0.56932 | 0.56932 | 0.56932 | 8.5 E−6 | 2.2 E−7 |
|  | 3.0 | 0.22782 | 0.22779 | 0.22772 | 0.22774 | 0.22774 | 0.22774 | 1.6 E−5 | 1.9 E−8 |

*Table 5: Comparison of $L_2$ and $L_\infty$ errors of case 2 for various time and viscosity at $\Delta t = 0.001$, $N=80$.*

| T | $\nu=0.005$ | | $\nu=0.002$ | | $\nu=0.0001$ | |
|---|---|---|---|---|---|---|
|  | $L_2$ | $L_\infty$ | $L_2$ | $L_\infty$ | $L_2$ | $L_\infty$ |
| 5  | 2.875E-09 | 5.419E-09 | 1.764E-08 | 5.931E-08 | 0.08022   | 0.29237   |
| 10 | 5.180E-10 | 9.913E-10 | 2.577E-09 | 3.938E-09 | 6.029E-07 | 1.246E-06 |
| 15 | 1.986E-10 | 3.762E-10 | 1.854E-09 | 3.296E-09 | 2.386E-07 | 5.058E-07 |

*Table 6: Absolute errors of numerical solutions with Δt=0.001, Re =20 and grid size of 16×16 at some specific points for 2D Burgers' model.*

| x | y | Liu et al. [53] | Proposed scheme | Liu et al. [53] | Proposed scheme | Liu et al. [53] | Proposed scheme |
|---|---|---|---|---|---|---|---|
| | | T=0.5 | | T=0.75 | | T=1.0 | |
| 0.125 | 0.125 | 1.56E−05 | 4.50E-06 | 1.21E−06 | 1.64E-06 | 3.93E−07 | 3.37E-07 |
| | 0.5 | 6.68E−05 | 4.92E-06 | 1.33E−05 | 3.40E-06 | 3.39E−06 | 8.56E-07 |
| | 0.875 | 1.01E−05 | 9.41E-07 | 1.60E−05 | 1.45E-06 | 1.87E−04 | 3.85E-06 |
| 0.5 | 0.125 | 6.68E−05 | 4.92E-06 | 1.33E−05 | 3.40E-06 | 3.39E−06 | 8.56E-07 |
| | 0.5 | 7.69E−07 | 5.60E-07 | 9.33E−06 | 2.01E-08 | 1.92E−05 | 6.11E-06 |
| | 0.875 | 7.91E−08 | 4.46E-08 | 1.12E−06 | 5.91E-07 | 5.82E−06 | 7.69E-07 |
| 0.875 | 0.125 | 1.01E−05 | 9.41E-07 | 1.60E−05 | 1.45E-06 | 1.87E−04 | 3.85E-06 |
| | 0.5 | 7.91E−08 | 4.46E-08 | 1.12E−06 | 5.91E-07 | 5.82E−06 | 7.69E-07 |
| | 0.875 | 3.47E−09 | 3.26E-09 | 4.80E−08 | 8.71E-09 | 5.23E−07 | 3.21E-07 |

*Table 7:  $L_2$ and $L_\infty$ error norms at T= 0.05, 0.25 and Re =1 for 2D Burgers' model.*

| Grid size | $\Delta t$ | Mittal and Tripathi [54] | | Liu et al. [53] | | Arora and Joshi [55] | | Proposed scheme | |
|---|---|---|---|---|---|---|---|---|---|
| | | $L_2$ | $L_\infty$ | $L_2$ | $L_\infty$ | $L_2$ | $L_\infty$ | $L_2$ | $L_\infty$ |
| T = 0.05 | | | | | | | | | |
| 5×5 | 0.005 | 4.969E−008 | 4.651E−008 | -- | 1.43E−06 | -- | -- | 4.375E-07 | 5.855E-07 |
| 10×10 | 0.0005 | 6.217E−009 | 5.907E−009 | -- | 2.40E−08 | -- | -- | 4.775E-09 | 4.492E-09 |
| 15×15 | 0.0001 | 2.532E−009 | 2.188E−009 | -- | 4.44E−09 | -- | -- | 2.407E-10 | 1.887E-10 |
| 30×30 | 0.0001 | 1.875E−009 | 1.010E−009 | -- | 5.58E−10 | -- | -- | -- | -- |
| 60×60 | 0.0001 | -- | -- | -- | -- | 5.62E-06 | 8.67E-06 | -- | -- |
| T = 0.25 | | | | | | | | | |
| 5×5 | 0.005 | 9.98993e−009 | 9.80769e−009 | -- | 9.11E−07 | -- | -- | 2.909E-07 | 4.057E-07 |
| 10×10 | 0.0005 | 8.12715e−009 | 7.04501e−009 | -- | 2.37E−08 | -- | -- | 2.379E-10 | 2.160E-10 |
| 15×15 | 0.0001 | 7.23104e−009 | 6.05481e−009 | -- | 4.56E−09 | -- | -- | 1.207E-11 | 8.888E-12 |
| 30×30 | 0.0001 | 3.93668e−009 | 3.06006e−009 | -- | 6.21E−10 | -- | -- | -- | -- |
| 60×60 | 0.0001 | -- | -- | -- | -- | 3.88E-06 | 4.99E-06 | -- | -- |

*Table 8: $L_2$ and $L_\infty$ error norms with $\Delta t=0.0005$ and grid size= 16×16 at different time and Re for 2D Burgers' model.*

| T | Re =10 | | Re =100 | | Re =200 | |
|---|---|---|---|---|---|---|
| | $L_2$ | $L_\infty$ | $L_2$ | $L_\infty$ | $L_2$ | $L_\infty$ |
| 3 | 3.18E−09 | 3.52E−09 | 3.11E−06 | 3.84E−06 | 1.56E−04 | 2.46E−04 |
| 5 | 1.29E−13 | 1.31E−13 | 1.35E−12 | 1.26E−12 | 1.93E−09 | 2.34E−09 |
| 10 | 3.08E−13 | 3.86E−13 | 7.59E−13 | 8.10E−13 | 9.52E−13 | 9.54E−13 |

*Table 9: $L_2$ and $L_\infty$ error norms for u component at Re =100 and different time.*

| T | Li et al.[56], N=64 $\Delta t$= 0.00390625 | | Proposed Meth., N=20 $\Delta t$= 0.001 | |
|---|---|---|---|---|
| | $L_2$ | $L_\infty$ | $L_2$ | $L_\infty$ |
| 0.5 | 7.4281E−04 | 18.698E−04 | 1.3078E−05 | 1.0721E−05 |
| 1 | 1.1101E−04 | 4.4963E−04 | 1.0779E−05 | 8.3286E−06 |
| 2 | 9.5229E−05 | 4.8575E−04 | 1.0823E−05 | 9.0187E−06 |
| 4 | 9.7265E−06 | 1.4358E−04 | 7.3885E−08 | 8.4375E−08 |

*Table 10: Comparisons with respect to exact solution and other Refs. of u at Re =100 and 20×20 grid size for case 4.3*

|  | T=0.5 | | | | T=2.0 | | | |
|---|---|---|---|---|---|---|---|---|
| (x , y) | Shukla et al. [57] Δt= 0.0001 | Shi et al. [58] | **Proposed Meth.** **Δt= 0.001** | Exact | Shukla et al. [57] Δt= 0.0001 | Shi et al. [58] | **Proposed Meth.** **Δt= 0.001** | Exact |
| (0.1,0.1) | 0.54412 | 0.54336 | 0.54332 | 0.54332 | 0.50050 | 0.50048 | 0.50048 | 0.50048 |
| (0.5,0.1) | 0.50037 | 0.50035 | 0.50035 | 0.50035 | 0.50000 | 0.50000 | 0.50000 | 0.50000 |
| (0.9,0.1) | 0.50000 | 0.50000 | 0.50000 | 0.50000 | 0.50000 | 0.50000 | 0.50000 | 0.50000 |
| (0.3,0.3) | 0.54388 | 0.54338 | 0.54338 | 0.54338 | 0.50050 | 0.50047 | 0.50048 | 0.50048 |
| (0.7,0.3) | 0.50037 | 0.50034 | 0.50035 | 0.50035 | 0.50000 | 0.50000 | 0.50000 | 0.50000 |
| (0.1,0.5) | 0.74196 | 0.74228 | 0.74222 | 0.74221 | 0.55632 | 0.55572 | 0.55568 | 0.55568 |
| (0.5,0.5) | 0.54347 | 0.54333 | 0.54332 | 0.54332 | 0.50050 | 0.50046 | 0.50048 | 0.50048 |
| (0.9,0.5) | 0.50035 | 0.50034 | 0.50035 | 0.50035 | 0.50001 | 0.50000 | 0.50000 | 0.50000 |
| (0.3,0.7) | 0.74211 | 0.74236 | 0.74223 | 0.74223 | 0.55597 | 0.55573 | 0.55577 | 0.55577 |
| (0.7,0.7) | 0.54327 | 0.54332 | 0.54338 | 0.54338 | 0.50054 | 0.50045 | 0.50048 | 0.50048 |
| (0.1,0.9) | 0.74994 | 0.74995 | 0.74995 | 0.74995 | 0.74406 | 0.74431 | 0.74426 | 0.74426 |
| (0.5,0.9) | 0.74219 | 0.74237 | 0.74221 | 0.74221 | 0.55575 | 0.55568 | 0.55568 | 0.55568 |
| (0.9,0.9) | 0.54333 | 0.54331 | 0.543324 | 0.543325 | 0.50052 | 0.50045 | 0.50048 | 0.50048 |

Table 11: Comparisons with respect to exact solution and other Refs. of v at Re =100 and 20×20 grid size for case 4.3

| (x , y) | T=0.5 | | | | T=2.0 | | | |
|---|---|---|---|---|---|---|---|---|
| | Shukla et al. [57] $\Delta t= 0.0001$ | Shi et al. [58] | Proposed Meth. $\Delta t= 0.001$ | Exact | Shukla et al. [57] $\Delta t= 0.0001$ | Shi et al. [58] | Proposed Meth. $\Delta t= 0.001$ | Exact |
| (0.1,0.1) | 0.95589 | 0.95665 | 0.95668 | 0.95668 | 0.99950 | 0.99952 | 0.99952 | 0.99952 |
| (0.5,0.1) | 0.99963 | 0.99965 | 0.99965 | 0.99965 | 1.00000 | 1.00000 | 1.00000 | 1.00000 |
| (0.9,0.1) | 1.00000 | 1.00000 | 1.00000 | 1.00000 | 1.00000 | 1.00000 | 1.00000 | 1.00000 |
| (0.3,0.3) | 0.95612 | 0.95662 | 0.95662 | 0.95662 | 0.99950 | 0.99953 | 0.99952 | 0.99952 |
| (0.7,0.3) | 0.99964 | 0.99966 | 0.99965 | 0.99965 | 1.00000 | 1.00000 | 1.00000 | 1.00000 |
| (0.1,0.5) | 0.75804 | 0.75772 | 0.75778 | 0.75779 | 0.94368 | 0.94428 | 0.94433 | 0.94432 |
| (0.5,0.5) | 0.95654 | 0.95667 | 0.95668 | 0.95668 | 0.99950 | 0.99954 | 0.99952 | 0.99952 |
| (0.9,0.5) | 0.99965 | 0.99966 | 0.99965 | 0.99965 | 0.99999 | 1.00000 | 1.00000 | 1.00000 |
| (0.3,0.7) | 0.75789 | 0.75764 | 0.75777 | 0.75777 | 0.94403 | 0.94427 | 0.94423 | 0.94423 |
| (0.7,0.7) | 0.95673 | 0.95668 | 0.95662 | 0.95662 | 0.99946 | 0.99955 | 0.99952 | 0.99952 |
| (0.1,0.9) | 0.75006 | 0.75005 | 0.75005 | 0.75005 | 0.75595 | 0.75569 | 0.75574 | 0.75574 |
| (0.5,0.9) | 0.75781 | 0.75763 | 0.75779 | 0.75779 | 0.94425 | 0.94432 | 0.94432 | 0.94432 |
| (0.9,0.9) | 0.95667 | 0.95669 | 0.956674 | 0.956675 | 0.99948 | 0.99955 | 0.99986 | 0.99986 |